\documentclass[11pt]{article}
\usepackage[latin1]{inputenc}
\usepackage[english]{babel}
\usepackage{amssymb,amsfonts, amsmath, color}
\usepackage[margin=2.7cm]{geometry}
\usepackage{xcolor}
\usepackage{graphicx, color, enumerate}
\usepackage[latin1]{inputenc}
\usepackage[active]{srcltx}
\usepackage{tikz}
\usepackage{pgf}
\usepackage{etex}
\usepackage{verbatim}
\usepackage{tikz-3dplot}
\usepackage{pgfkeys}
\usepackage{amsmath}
\usepackage{amsfonts}
\usepackage{amssymb}
\usepackage{amsthm}
\usepackage{algpseudocode}
\usepackage{float}
\usepackage{xcolor}
\usepackage{mathrsfs}
\usepackage{import}
\usepackage{geometry}
\usepackage{fancyhdr}
\usepackage{fp}

\newtheorem{theorem}{Theorem}[section]
\newtheorem{proposition}[theorem]{Proposition}
\newtheorem{corollary}[theorem]{Corollary}
\newtheorem{lemma}[theorem]{Lemma}
\newtheorem{remark}[theorem]{Remark}
\newtheorem{definition}[theorem]{Definition}

\newcommand{\bcl}{\begin{center}}
\newcommand{\ecl}{\end{center}}
\newcommand{\brl}{\begin{right}}
\newcommand{\erl}{\end{right}}
\newcommand{\ben}{\begin{enumerate}}
\newcommand{\een}{\end{enumerate}}
\newcommand{\overliner}{\begin{array}}
\newcommand{\earr}{\end{array}}
\newcommand{\btab}{\begin{tabular}}
\newcommand{\etab}{\end{tabular}}
\newcommand{\bdoc}{\begin{document}}
\newcommand{\edoc}{\end{document}}
\newcommand{\beqy}{\begin{eqnarray}}
\newcommand{\eeqy}{\end{eqnarray}}

\newcommand{\beqi}{\begin{eqnarray*}}
\newcommand{\eeqi}{\end{eqnarray*}}
\newcommand{\bitem}{\begin{itemize}}

\newcommand{\eitem}{\end{itemize}}
\newcommand{\nln}{\newline}
\newcommand{\newt}{\newtheorem}

\newcommand{\pa}{\partial}
\newcommand{\re}{{I\!\!R}}
\newcommand{\Rn}{\R^N}
\newcommand{\xr}{x\in\R }
\newcommand{\x}{\times}
\newcommand{\dyle}{\displaystyle}
\newcommand{\ene}{{I\!\!N}}
\newcommand{\irn}{\int\limits_{\R^N}}
\newcommand{\io}{\int\limits_{\O}}
\newcommand{\meas}{{\rm meas\,}}

\newcommand{\sign}{{\rm sign}}
\newcommand{\map}{\longrightarrow }
\newcommand{\imp}{\Longrightarrow }
\renewcommand{\div}{\nabla\cdot }
\newcommand{\sen}{{\rm sen\,}}
\newcommand{\tg}{{\rm tg\,}}
\newcommand{\arcsen}{{\rm arcsen\,}}
\newcommand{\arctg}{{\rm arctg\,}}
\newcommand{\supp}{{\textsl supp\ }}
\newcommand{\ity}{\int_{-\iy}^{+\iy}}
\newcommand{\limit}{\lim\limits}
\newcommand{\limi}{\limit_{n\to\infty}}
\newcommand{\sumi}{\sum\limits_{n=1}^{\infty}}
\newcommand{\ulu}{\underline u}
\newcommand{\ulw}{\underline w}
\newcommand{\ulz}{\underline z}
\newcommand{\ulv}{\underline v}
\newcommand{\uls}{\underline s}
\newcommand{\olu}{\overline u}
\newcommand{\olv}{\overline v}
\newcommand{\ols}{\overline s}
\newcommand{\ob}{\overline\b}
\newcommand{\ovar}{\overline\var}
\newcommand{\wv}{\widetilde v}
\newcommand{\wu}{\widetilde u}
\newcommand{\ws}{\widetilde s}
\renewcommand{\a }{\alpha }
\renewcommand{\b }{\beta }
\newcommand{\g }{\gamma}
\newcommand{\G }{\Gamma }
\renewcommand{\d }{\delta }

\newcommand{\D }{\Delta }
\newcommand{\e }{\varepsilon }
\newcommand{\z }{\zeta }
\renewcommand{\l }{\lambda }
\renewcommand{\L }{\Lambda }
\newcommand{\m }{\mu }
\newcommand{\n }{\nabla }
\newcommand{\s }{\sigma }
\newcommand{\Sig }{\Sigma }
\renewcommand{\t }{\tau }
\newcommand{\var }{\varphi }
\renewcommand{\o }{\omega }
\renewcommand{\O }{\Omega }
\newcommand{\R}{{\mathbb{R}}}
\newcommand{\bC}{{\bf C}}
\newcommand{\bZ}{{\bf Z}}
\newcommand{\bN}{{\bf N}}
\newcommand{\bQ}{{\bf Q}}
\newcommand{\bK}{{\bf K}}
\newcommand{\bI}{{\bf I}}
\newcommand{\bv}{{\bf v}}
\newcommand{\bV}{{\bf V}}
\DeclareMathOperator{\suppo}{supp} \DeclareMathOperator{\di}{div}

\newenvironment{Proof}{\Rmovelastskip\vskip12pt
plus 1pt \noindent\em\rm}{\hfill {\qed \hskip .2cm}}

\title{Uniqueness for fractional parabolic and elliptic equations with drift}
\author{Giulia Meglioli \thanks{Dipartimento di Matematica, Politecnico di Milano, (giulia.meglioli@polimi.it).}\; 
and Fabio Punzo\thanks{Dipartimento di Matematica, Politecnico di Milano, (fabio.punzo@polimi.it).}\; }

\date{}

\begin{document}
\maketitle

\abstract{We investigate uniqueness, in suitable weighted Lebesgue spaces, of solutions to a class of fractional parabolic and elliptic equations with a drift.}

\bigskip
\smallskip



\section{Introduction}\setcounter{equation}{0}
We are concerned with uniqueness of solutions, in suitable weighted Lebesgue spaces, to the following linear nonlocal Cauchy problem:
\begin{equation}\label{problema}
\begin{cases}
\rho\, u_t+(-\Delta)^s u - \left \langle b,\,\nabla u \right \rangle =0 \quad & \text{in}\,\,\,S_T:=\R^N\times(0,T]\\
u=0\quad & \text{in}\,\,\,\R^N\times\{0\}.
\end{cases}
\end{equation}
where $\rho$, usually referred to as a variable density, is a positive function depending only on the space variable, $(-\Delta)^s$ denotes the fractional Laplace operator of order $s\in (0,1)$ and $b:\R^N\to\R^N$ is a given vector field.
Moreover, we investigate existence and uniqueness of solutions to the linear nonlocal elliptic equation
\begin{equation}\label{elliptic}
(-\Delta)^{s} u - \left \langle b,\,\nabla u \right \rangle + \rho c  u \,=\,0 \quad \textrm{in}\;\; \R^N,
\end{equation}
where $c$ is a nonnegative function defined in $\R^N$. 

\medskip

We always assume that 
\begin{equation}\tag{{\it $H_0$}}\label{h3}
\begin{aligned}
&\textrm{(i)} \,\rho\in C(\R^N);\\
&\textrm{(ii)} \;\text{there exist}\,\,\alpha\ge 0\,\,\text{and}\,\,C_0>0\,\,\text{such that}\,\,\\
&\quad \; \rho(x)\geq C_0 (1+|x|^2)^{-\frac{\alpha}2} \quad\text{for all}\,\,x\in \mathbb R^N.
\end{aligned}
\end{equation}

Set
\begin{equation}\label{d}
\begin{aligned}
&D_+:=\{x\in\R^N\,\,:\, \langle b(x),x\rangle> 0\};\\
&D_-:=\{x\in\R^N\,\,:\, \langle b(x),x\rangle\le 0\}.
\end{aligned}
\end{equation}
Concerning the vector field $b$, we always make the following assumptions:
\begin{equation}\tag{{\it $H_1$}}\label{h0}
\begin{aligned}
&\textrm{(i)} \; b:\R^N\to\R^N,\,\,\,b\in C^1(\R^N)\,;\\
&\textrm{(ii)} \;\text{there exist}\,\,\sigma\le 1-\alpha\,\,\text{and}\,\,K>0\,\,\text{such that}\,\,\\
&\quad\quad \quad\left\langle b(x), \frac{x}{|x|}\right\rangle\le K(1+|x|)^{\sigma}\quad\text{for all}\,\,x\in D_+;\\
&\textrm{(iii)}\;[\operatorname{div} b(x)]_{-}\le K(1+|x|)^{\sigma-1}\quad\text{for all}\,\,x\in\R^N;
\end{aligned}
\end{equation}
Here $\alpha$ is the same as in \eqref{h3}; in addition, for a given function $v$, the negative part $[v]_{-}$ is 
$$
[v]_{-}:=\max\left\{0;\,-v\right\}.
$$
Moreover, the coefficient $c$ is such that
\begin{equation}\label{h1}\tag{{\it $H_2$}}
c\in C(\R^N),\; c(x)\geq 0 \quad \textrm{for all}\;\; x\in \R^N\,.
\end{equation}

The fractional diffusion equation with $\rho\equiv 1$ and a drift term $b$, appearing in \eqref{problema}, has been the object of an important investigation in the literature. In particular, in \cite{CaffVass}, \cite{Silv}, \cite{Silv2} and in \cite{Silv3} regularity of solutions is addressed. Moreover, in \cite{CW}, \cite{Zhang1}, \cite{Zhang2} and in \cite{CSZ} the relation between problem \eqref{problema} and L\'evy processes has been studied.

Furthermore, in the literature there are many results about uniqueness of solutions to certain problems somehow related to \eqref{problema} and \eqref{elliptic}.
We now briefly recall some of such known results. In doing it, we distinguish two cases: (a) the local version of  \eqref{problema} and \eqref{elliptic}; (b) nonlocal equations.

\medskip

\noindent{\it (a) The local case.} In \cite{AB}, the operator
$$
L u:= \sum_{i,j=1}^N\frac{\pa^2 \big[a_{ij}(x,t) u\big] }{\pa x_i\pa x_j} - \sum_{i=1}^N\frac{\pa \big[b_i(x,t) u\big]}{\pa x_i}+ c(x,t)u - u_t\,
$$
is considered; the coefficients of $L$, together with all their derivatives which appear, are locally bounded functions in $S_T$.
Furthermore, the matrix $A\equiv (a_{ij})$ is assumed to be positive semidefinite in $\bar S_T$.  Suppose that
$$
|a_{ij}(x,t)|\leq K_1 ( 1+|x|^2)^{\frac{2-\l}{2}},\; |b_i(x,t)|\leq K_2(|x|^2+1)^{\frac{1}{2}},\;|c(x,t)|\leq K_3(|x|^2+1)^{\frac{\l}2},
$$
for almost every $(x,t)\in S_T$, for some constants $\l\geq 0, K_i>0\; (i=1,2,3)$.
\smallskip

\noindent For any $\phi\in C(\overline S_T)$, $\phi>0$, set
$$
L^1_\phi(S_T):= \left\{u: S_T\to \R\,\;\textrm{measurable}\;\,|\, \int_0^T\int_{\R^N}|u|\phi(x,t) dx dt<\infty\right\}\,.
$$
In \cite[Theorem 1]{AB} it is shown that if $u$ is a solution to problem 
 \begin{equation}\label{eq18}
\begin{cases}
L\, u = 0\quad &\textrm{in}\,\,S_T\\
u \, = 0 \quad&\textrm{in\ \ } \R^N\times \{0\} \,,
\end{cases}
\end{equation}
and $u\in L^1_{\phi}(S_T)$, with
\begin{equation}\label{eq19}
\phi(x)=(|x|^2+1)^{-\a_0} \quad (x\in \R^N)\;\; \textrm{if}\;\;
\l=0,
\end{equation}
or
\begin{equation}\label{eq15}
\phi(x)=e^{-\a_0(|x|^2+1)^{\frac{2-\a}2}}\quad (x\in \R^N)\;\;
\textrm{if}\;\; \l>0\,,
\end{equation}
for some $\a_0>0$, then
$$
u\equiv 0\quad \textrm{in}\;\; S_T\,.
$$


Clearly, as a consequence of the previous result one can immediately deduce uniqueness in $L^1_{\phi}(S_T)$ for solutions to problem
$$
\begin{cases}
L u = f\quad &\textrm{in}\,\,S_T\\
u = u_0 \quad&\textrm{in\ \ } \R^N\times \{0\} \,,
\end{cases}
$$
where $f$ and $u_0$ are given functions defined in $S_T$ and $\R^N$, respectively.

\bigskip

Furthermore, in the literature (see e.g. \cite{AB, EKP, IKO, KPT}),  it has been widely studied uniqueness of solutions to the following problem
\begin{equation*}
\begin{cases}
\rho u_t -\Delta u  =0 \quad & \text{in}\,\,\,S_T\\
u=0\quad & \text{in}\,\,\,\R^N\times\{0\}.
\end{cases}
\end{equation*}

Concerning the elliptic equation \eqref{elliptic} with $s=1$, some existence, uniqueness and nonunqueness results have been established in \cite{OR}, \cite{PTes} and in \cite{P2}. In particular,  it is addressed the case of the equation posed in bounded domains, with coefficients that can be degenerate or singular at the boundary of the domain.

\bigskip

\noindent {\it (b) The fractional case.} In \cite{PV1} and in \cite{PV2}, it is studied problem
\begin{equation}\label{eq120}
\begin{cases}
\rho u_t+(-\Delta)^s u  =0 \quad & \text{in}\,\,\,S_T\\
u=0\quad & \text{in}\,\,\,\R^N\times\{0\}\,,
\end{cases}
\end{equation}
that is problem \eqref{problema} with $b\equiv 0$. In particular, it has been investigated how $\rho$ affects uniqueness and nonuniqueness of  solutions. Uniqueness is obtained for solutions belonging to $L^p_\psi(S_T)$, for suitable weight $\psi$, depending on the behaviour of $\rho$ as $|x|\to +\infty$. 

Similar results have been also obtained for equation
$$
(-\Delta)^{s} u + \rho \,c u \,=\,0 \quad \textrm{in}\;\; \R^N\,.
$$

Finally, let us mention that in \cite[Lemma 2.4]{CSZ}, by means of probabilistic methods, it is shown that problem \eqref{problema} with $\rho\equiv 1$ admits a unique {\it bounded} solution, provided that $b$ is {\it bounded}.

\subsection{Outline of our results} 
The main results of this paper will be given in detail in the forthcoming Theorems \ref{teo1} and \ref{teo3}. We give here a sketchy outline of these results. 

\medskip

We shall prove that the solution to problem \eqref{problema} is unique in the class $L^p_{\psi}(S_T)$ with $p\geq 1$ and
\begin{equation}\label{eq12}
\psi(x):=(1+|x|^2)^{-\frac{\beta}2}\;\,\; (x\in \R^N)\,,
\end{equation}
for properly chosen $\beta>0$, provided that $b$ satisfies assumption \eqref{h0}. Such hypothesis allows that $b$ is {\it unbounded}. In addition, a growth condition on $b$ is required, when $\langle b(x), x\rangle> 0$. However, when $\langle b(x), x\rangle\leq 0$ no further conditions on $b$ are imposed.

In order to prove such a uniqueness result we construct a positive supersolution to equation
\begin{equation}\label{eq14}
- \rho\, \phi_t -(-\Delta)^{s} \phi- \left \langle b,\,\nabla u \right \rangle - \phi\di b= 0\quad \textrm{in}\;\; S_T\,.
\end{equation}
Indeed, the weight function $\psi$ defined in \eqref{eq12} is related to such a supersolution. In general, our uniqueness class includes {\it unbounded} solutions. Thus, in particular, we get uniqueness of bounded solutions.  Furthermore, we show sharpness of the hypothesis on the drift term $b$ (see Proposition \ref{cor2}). More precisely, we show that:
\begin{itemize}
\item if $b:\R^N\to\R^N$ satisfies assumption \eqref{h0}, then uniqueness of the solution to problem \eqref{problema} is guaranteed in $L^p_{\psi}$ for $p\ge 1$, for a suitable weight $\psi$ (see Theorem \ref{teo1} below);

\item conversely, if the drift term $b$ violetes, in an appropriate sense, assumption \eqref{h0} (that is, if it satisfies inequality \eqref{h2} below), then infinitely many bounded solutions to problem \eqref{problema} exist (see Proposition \ref{cor2} below).
\end{itemize}

Clearly, from our uniqueness result we can infer also uniqueness in $L^p_{\psi}(S_T)$ of solutions to
\begin{equation}\label{eq13}
\begin{cases}
\rho u_t+(-\Delta)^s u - \left \langle b,\,\nabla u \right \rangle =f \quad & \text{in}\,\,\,S_T\\
u=u_0\quad & \text{in}\,\,\,\R^N\times\{0\}.
\end{cases}
\end{equation}

\smallskip

Finally, we also establish similar uniqueness results for the linear elliptic {\it nonlocal} equation \eqref{elliptic} (see Theorem \ref{teo3} below). 
\medskip

Observe that our results generalize those concerning uniqueness established in \cite{CSZ}, recalled above. In fact, we consider possible unbounded solutions, an unbounded drift term $b$ and a density $\rho$. Let us underline that our methods are completely different from those used in \cite{CSZ}. 

\medskip

The paper is organized as follows. In Section \ref{sec1} we recall some preliminaries about fractional Laplacian and we give the notion of solutions we shall deal with. Then we state our main results concerning both parabolic and elliptic problems. Section \ref{pp} is devoted to the proof of results for parabolic problems, instead those about elliptic equations are proved in Section \ref{pe}.

\section{Mathematical framework and results} \label{sec1}\setcounter{equation}{0}
The fractional Laplacian $(-\Delta)^{s}$ can be defined by Fourier transform. Namely, for any function $g$ in the Schwartz class
$\mathcal S$, we say that
$$(-\Delta)^{s/2} g = h\,,$$
if
\begin{equation}
\label{eq21}
\hat{h}(\xi)= |\xi|^{s}\hat{g}(\xi).
\end{equation}
Here, we used the notation~$\hat{h}={\mathfrak F} h$ for the Fourier transform of~$h$. Furthermore, consider the space
$$
\mathcal L^s(\R^N):= \left\{u:\R^N\to \R\,\;\textrm{measurable}\, \,:\,\int_{\R^N}\frac{|u(x)|}{1+|x|^{N+2s}}dx<\infty \right\}\,, 
$$
 endowed with the norm
$$
 \| u \|_{\mathcal L^s(\R^N)}:=\int_{\R^N}\frac{|u(x)|}{1+|x|^{N+2s}}dx\,.
 $$
If $u\in \mathcal L^s(\R^N)$ (see \cite{Silv}), then $(-\Delta)^s u$ can be defined as a distribution, $i.e.$, for any $\varphi\in\mathcal S$,
$$
 \int_{\R^N} \varphi (-\Delta)^s u \, dx \,=\,   \int_{\R^N}  u (-\Delta)^s \varphi \, dx\,. 
$$
In addition, suppose that, for some $\g>0$, $u\in \mathcal L^s(\R^N)\cap C^{2s+\g}(\R^N)$ if $s<\frac 1 2$, or  $u\in \mathcal L^s(\R^N)\cap C^{1,
2s+\g-1}_{loc}(\R^N)\; \textrm{if}\; s\geq \frac 1 2$. Then we have
\begin{equation}\label{eq22}
(-\Delta)^{s} u(x)=C_{N,s}\,\, \textrm{P.V.}\, \int_{\R^N} \frac{u(x)-u(y)}{|x-y|^{N+2s}}d y\quad (x\in \R^N),
\end{equation}
where 
$$
C_{N,s}=\frac{2^{2s-1}{2s}\Gamma((N+2s)/2)}{\pi^{N/2}\Gamma(1-s)},
$$ 
$\Gamma$ being the Gamma function; moreover, $(-\Delta)^s u\in C(\R^N)$. Note that, see \cite{DPVal}, the constant $C_{N,s}$ satisfies the identity
\begin{equation}\label{eq23}
(-\Delta)^s u = \mathfrak F^{-1} \big( |\xi|^{2s}\mathfrak F u  \big)\,,\quad \xi\in\R^N, u\in \mathcal S\,,
\end{equation}
hence
$$
C_{N,s}=\left(\int_{\R^N}\frac{1-\cos(\xi_1)}{|\xi|^{N+2s}}d\xi \right)^{-1}\,.
$$

Now we can give the definition of solution to problem \eqref{problema} and to equation \eqref{elliptic}.

\begin{definition}\label{defsoleqp}
We say that a function $u$ is a {\em solution} to equation
\begin{equation}\label{eq24}
\rho\, u_t+(-\Delta)^s u - \left \langle b,\,\nabla u \right \rangle =\,0\quad \textrm{in}\;\; S_T\,,
\end{equation}
if
\begin{itemize}
\item[(i)]  $u\in C(S_T)$, for each $t\in (0,T]$  $u(\cdot, t)\in \mathcal L^s(\R^N)\cap
C^{2s+\g}_{loc}(\R^N)$ if $s<\frac 1 2$, or  $u(\cdot, t)\in \mathcal L^s(\R^N)\cap C^{1,2s+\g-1}_{loc}(\R^N)\; \textrm{if}\; s\geq \frac 1 2$,  for some $\g>0$, $\pa_t u\in C(S_T)$\,;
\item[(ii)] $ \rho(x) u_t + C_{N,s} \textrm{P.V}.\, \displaystyle\int_{\R^N} \frac{u(x,t)-u(y,t)}{|x-y|^{N+2s}}d y \,- \left \langle b,\,\nabla u \right \rangle=\,0$\;\;  for all\;\; $(x,t)\in
S_T$\,.
\end{itemize}
Furthermore, we say that $u$ is a {\em supersolution\; (subsolution)} to equation \eqref{eq24}, if in $(ii)$ instead of $``="$ we have $``\geq"\; (``\leq")$\,.
\end{definition}

\begin{definition}\label{defsolp}
We say that a function $u$ is a solution to problem \eqref{problema} if
\begin{itemize}
\item[(i)] $u\in C(\bar S_T), u_t\in L^1_{loc}(\bar S_T), u\in
L ^1\big((0,T),\mathcal L^s(\R^N)\big)$\,;
\item[(ii)]  $u$ is a solution to equation \eqref{eq24}\, in the sense of Definition \ref{defsoleqp};
\item[(iii)] $u(x,0)=0$\;\; for all\;\; $x\in \R^N\,.$
\end{itemize}
\end{definition}

\begin{definition}\label{defsole}
We say that a function $u$ is a solution to equation \eqref{elliptic} if
\begin{itemize}
\item[(i)]  $u\in \mathcal L^s(\R^N)\cap C^{2s+\g}_{loc}(\R^N)$ if $s<\frac 1 2$, or  $u\in \mathcal L^s(\R^N)\cap C^{1,2s+\g-1}_{loc}(\R^N)\; \textrm{if}\; s\geq \frac 1 2$, for some $\g>0$\,;
\item[(ii)]  $C_{N,s} \textrm{P.V}.\, \displaystyle\int_{\R^N} \frac{u(x)-u(y)}{|x-y|^{N+2s}}d y - \left \langle b,\,\nabla u \right \rangle + \rho(x) c(x)u(x) \,=\,0$\;\;  for all\;\; $x\in \R^N$\,.
\end{itemize}
Furthermore, we say that $u$ is a {\em supersolution\; (subsolution)} to equation \eqref{elliptic}, if in $(ii)$ instead of $``="$ we have $``\geq"\; (``\leq")$\,.
\end{definition}

\subsection{Parabolic equations: results} 

\begin{theorem}\label{teo1}
Let assumptions \eqref{h3} and \eqref{h0} be satisfied. Let $u$ be a solution to problem \eqref{problema} with $|u(\cdot, t)|^p\in \mathcal L^s(\R^N)$, for some
$p\geq 1$, for each $t>0$. Assume that one of the following condition is fulfilled:
\begin{itemize}
\item[(i)] $0<\beta\leq N-2s$, $\alpha\geq 0$;
\item[(ii)] $N-2s<\beta<N$, $\alpha\leq 2s$;
\item[(iii)] $\beta=N$, $\alpha<2s$;
\item[(iv)] $\beta>N$, $\alpha+\beta\leq 2s +N$.
\end{itemize} 
Let $\psi$ be defined as in \eqref{eq12}. If $u\in L^p_{\psi}(S_T)$, then
$$
u\equiv 0 \quad \textrm{in} \;\; S_T.
$$
\end{theorem}
\medskip


From Theorem \ref{teo1} we deduce the following
\medskip

\begin{corollary}\label{cor1}
Let assumption $(H_0)$ be satisfied. Let $u$ be a solution to problem \eqref{problema}. Suppose that $\alpha<2s$. If 
$$|u(x,t)| \leq C(1+|x|^2)^{\frac{\alpha}2}\quad \textrm{for all}\;\; x\in S_T,$$ 
for some $\alpha\in (0, 2s)$ and $C>0$, then
$$u\equiv 0\quad \textrm{in}\;\; S_T\,.$$
\end{corollary}

In order to prove Corollary \ref{cor1} it suffices to apply Theorem \ref{teo1} with $\beta=N+2s-\alpha>N$ and $p=1$. 

\smallskip

\noindent Note that hypothesis \eqref{h0} in Theorem \ref{teo1} is optimal. In fact, in the next proposition, by choosing a vector field $b$ for which \eqref{h0} fails, we see that problem \eqref{problema} admits infinitely many bounded solutions.

\normalcolor 

\medskip
\begin{proposition}\label{cor2} Let $\rho\in C^1(\mathbb R^N)$, $\rho>0$ in $\R^N$.
Let $b\in C^1(\R^N)$ be such that
\begin{equation}\label{h2}
\begin{aligned}
\langle b(x), x \rangle \geq 0 \quad &\text{ for all }\,\, x\in \mathbb R^N\,,\\
\left\langle b(x),\frac{x}{|x|}\right\rangle \ge\, K|x|^{\sigma}\quad&\text{for all}\,\,x\in\R^N\setminus B_{R_0},
\end{aligned}
\end{equation}
for some $$R_0>0, \quad \sigma>1-\alpha \quad \text{and}\quad K>0.$$
Then problem \eqref{problema} admits infinitely bounded many solutions. 
\end{proposition}

\subsection{Elliptic equations: results}

\begin{theorem}\label{teo3}
Let assumptions \eqref{h3}-\eqref{h1} be satisfied. Let $u$ be a solution to equation \eqref{elliptic} with $|u|^p\in \mathcal L^s(\R^N),$ for some $p\geq 1$. Suppose that, for some $c_0>0,$
\begin{equation}\label{eq27}
c(x)\geq c_0\quad \textrm{for all}\;\;x\in \R^N\,.
\end{equation}
Assume that one of the conditions $(i)-(iv)$ of Theorem \ref{teo1} holds. Furthermore, suppose that  
that $p c_0$ is large enough. Let $\psi$ be defined as in \eqref{eq12}. If $u\in L^p_{\psi}(\R^N)$, then
$$
u\equiv 0 \quad \textrm{in} \;\; \R^N\,.
$$
\end{theorem}

\begin{remark}
{\rm The hypothesis $p c_0 $ large enough made in Theorem \ref{teo3} will be specified in the proof of Theorem \ref{teo3}.}
\end{remark}

Analogously to Corollary \ref{cor1}, we have the following
\begin{corollary}\label{cor3}
Let assumptions \eqref{h0} and \eqref{h1} be satisfied. Let $u$ be a solution to equation \eqref{elliptic} with $|u|^p\in \mathcal L^s(\R^N),$ for some $p\geq 1$ . Suppose that $pc_0$ is large enough, where $c_0$ has been defined in \eqref{eq27}. Suppose that $\alpha<2s$. If 
$$|u(x)| \leq C(1+|x|^2)^{\frac{\alpha}2}\quad \textrm{for all}\;\; x\in \R^N,$$ 
for some $\alpha\in (0, 2s)$ and $C>0$, then
$$u\equiv 0\quad \textrm{in}\;\; \R^N\,.$$
\end{corollary}

In order to prove Corollary \ref{cor3} it suffices to apply Theorem \ref{teo3} with $\beta=N+2s-\alpha>N$ and $p=1$.

\smallskip

Also for the elliptic equation \eqref{elliptic}, the hypothesis \eqref{h0} on $b$ is sharp. This is the content of the next proposition. 

\medskip
\begin{proposition}\label{cor4} Let $\rho\in C^1(\mathbb R^N)$, $\rho>0$ in $\R^N$.
Let  assumptions \eqref{h1} and \eqref{h2} be in force. Then equation \eqref{elliptic} admits infinitely many bounded solutions. 
\end{proposition}

\section{Parabolic equations: proofs}\setcounter{equation}{0} \label{pp}

Let us observe that 
\begin{equation*}
\begin{aligned}
&\text{if}\,\, f,g  \in \mathcal L^s(\R^N)\cap C_{loc}^{2s+\gamma}(\R^N)\,\,\,\,\quad\quad \text{with}\,\,\, s<\frac 1 2, \\
&\text{or}\,\,f,g\in \mathcal L^s(\R^N)\cap C^{1,2s+\gamma-1}_{loc}(\R^N)\,\,\,\, \text{with}\,\,\, s\geq \frac 1 2, 
\end{aligned}
\end{equation*}
for some $\gamma>0$, and $fg\in \mathcal L^s(\R^N)$, then it is easily checked that
\begin{equation}\label{eq30}
(-\Delta)^s[f(x)g(x)]= f(x)(-\Delta)^{s}g(x) + g(x)(-\Delta)^sf(x)- \mathcal B(f,g)(x),
\end{equation}
for all $x\in \R^N\,$, where $\mathcal B(f,g)$ is the bilinear form given by
$$
\mathcal B(f,g)(x):=C_{N,s}\int_{\R^N}\frac{[f(x)-f(y)][g(x)-g(y)]}{|x-y|^{N+2s}}dy \quad \textrm{for all}\;\; x\in \R^N\,.
$$
Take a cut-off function $\gamma\in C^\infty([0,\infty)), 0\leq \g\leq1$ with
\begin{equation}\label{eq31}
\gamma(r)=
\begin{cases}
1 &\textrm{if}\,\,0\leq r\leq \frac 1 2\\
0& \textrm{if\ \ } r\geq 1 \,
\end{cases};
\quad\quad\quad \gamma'(r)<0.
\end{equation}
Moreover, for any $R>0$ let
\begin{equation}\label{eq32}
\gamma_R(x):= \gamma\left(\frac{|x|}{R}\right) \quad \textrm{for all}\;\; x\in \R^N\,;
\end{equation}
and for any $\tau\in (0,T)$ let
$$
S_\tau:= \R^N\times (0,\tau]\,.
$$

\medskip

Next, we prove a general criterion for uniqueness of nonnegative solutions to problem
\eqref{problema} in $L^1_{\psi}(S_T),$ where $\psi$ is defined as in \eqref{eq12} for some constant $\beta>0$.

\begin{proposition}\label{prop1}
Let assumptions \eqref{h3}, \eqref{h0}-(i) be satisfied. Let $u$ be a solution to problem \eqref{problema} with $|u(\cdot, t)|^p\in \mathcal L^s(\R^N)$ for some $p\geq 1$, for each $t>0$. Assume that there exists a positive supersolution $\phi\in C^2(\bar S_T)$ to equation
\begin{equation}\label{eq32b}
\rho \phi_t-(-\Delta)^s\phi  - \left \langle b,\,\nabla \phi \right \rangle -\phi\di b =0 \quad \textrm{in}\;\; S_T\,,
\end{equation}
such that
\begin{equation}\label{eq32d}
\phi(x,t)+|\nabla \phi(x,t)|\leq C \psi(x)\quad \textrm{for all}\;\; (x,t)\in S_T,
\end{equation}
and
\begin{equation}\label{eq32c}
\frac{1}{1+|x|}\left\langle b(x),\frac{x}{|x|}\right\rangle\phi(x,t)\leq C \psi(x)\quad \textrm{for all}\;\; (x,t)\in D_+\times(0,T);
\end{equation}
for some constant $C>0$. If $u\in L^p_\psi(S_T)$, then
$$u\equiv 0\quad \textrm{in}\;\; S_T\,.$$
\end{proposition}

\subsection{Proof of Proposition \ref{prop1}}

To prove Proposition \ref{prop1} we shall use the next two results.

\begin{lemma}\label{lemma1}
Let $\tau\in (0,T)$, $\phi\in C^2(\bar S_\tau)$, $\phi>0$; suppose that \eqref{h0}-(i), \eqref{eq32d} and \eqref{eq32c} are satisfied. Let $v\in L^1_{\psi}(S_\t)$. Then
\begin{equation}\label{eq31bis}
\int_0^\t\int_{\R^N}|v(x,t)|\phi(x,t)|(-\Delta)^s\g_R(x)|\,dxdt +\int_0^\t\int_{\R^N}|v(x,t)|\,|\mathcal B(\phi, \g_R)(x)|dx dt\longrightarrow 0 \end{equation}
as $R\to\infty$, and
\begin{equation}\label{eq31tris}
\lim_{R\to\infty}\, \int_0^\t\int_{D_+}|v(x,t)|\phi(x,t) \left\langle b(x),\nabla\gamma_R(x)\right\rangle\,dx dt=0,
\end{equation}
where $D_+$ is defined in \eqref{d}.
\end{lemma}
\begin{remark}\label{rem21}
Observe that the quantity into the brackets of formula \eqref{eq31tris} is negative. In fact, due to \eqref{d} and \eqref{eq31}, for any $x\in D_+$,
$$
\left\langle b(x),\nabla\gamma_R(x)\right\rangle=\gamma'\left(\frac{|x|}{R}\right)\frac{1}{R}\,\left\langle b(x),\frac{x}{|x|}\right\rangle\le 0;
$$
and, by assumption in Lemma \ref{lemma1}, $\phi>0$ in $S_T$.
\end{remark}
Observe that a similar result was obtained in \cite{PV1}. However, in \cite{PV1}, only assumption \eqref{eq31bis} was made. This difference arises because of the presence of the vector field $b$ in \eqref{problema} which has not been considered in \cite{PV1}. For this reason, the proof of \ref{lemma1} is a slight modification of the proof of \cite[Lemma 3.1]{PV1}, hence we just show how to treat the extra term given by the vector field $b$.

\medskip

\begin{proof}[Proof of Lemma \ref{lemma1}]
In view of \eqref{eq32d}, by arguing as in the proof of \cite[Lemma 3.1]{PV1}, we get
\begin{equation}\label{eq33a}
\begin{aligned}
&\int_0^\t\int_{\R^N}|v(x,t)|\phi(x,t)|(-\Delta)^s\g_R(x)|\,dx dt \longrightarrow 0\\
&\int_0^\t\int_{\R^N}|v(x,t)|\,|\mathcal B(\phi, \g_R)(x)|\,dx dt \longrightarrow 0,
\end{aligned}
\end{equation}
as $R\to \infty$, hence \eqref{eq31bis} is proved. 

\noindent To show \eqref{eq31tris}, let us observe that, due to \eqref{eq31}, for any $x\in D_+\cap (B_R\setminus B_{R/2})$, for some $\bar C>0$
\begin{equation}\label{eq33b}
-\left\langle b(x),\nabla \gamma_R(x)\right\rangle= -\gamma'\left(\frac{|x|}{R}\right)\frac{1}{R}\,\left\langle b(x),\frac{x}{|x|}\right\rangle\,\le \bar C\frac{1}{|x|}\,\left\langle b(x),\frac{x}{|x|}\right\rangle,
\end{equation}
with $D_+$ as in \eqref{d}. Then, due to \eqref{eq32c} together with \eqref{eq33b}, we get
\begin{equation}\label{eq33}
\begin{aligned}
-\int_0^\tau\int_{D_+} &|v(x,t)|\phi(x,t)\langle b(x),\nabla \gamma_R\rangle\,dxdt  \\
&\le \bar C \int_0^{\tau}\int_{D_+\cap (B_R\setminus B_{R/2})} |v(x,t)|\,\phi(x,t)\,\frac{1}{|x|}\,\left\langle b(x),\frac{x}{|x|}\right\rangle\,dxdt \\
&\le \bar C C \int_0^{\tau}\int_{D_+\cap (B_R\setminus B_{R/2})} |v(x,t)| \psi(x)\,dxdt.
\end{aligned}
\end{equation}
Then, since $u\in L^1_{\psi}(S_T)$, we obtain from \eqref{eq33} that
\begin{equation}\label{eq34}
\lim_{R\to\infty}\int_0^{\tau}\int_{D_+} |v(x,t)|\,\phi(x,t)\,\left\langle b(x),\nabla \gamma_R(x)\right\rangle\,dxdt=0.
\end{equation}
This completes the proof.
\end{proof}

\begin{lemma}\label{lemma2}
Let $G\in C^2(\R; \R)$ be a convex function. Let $u\in \mathcal L^s(\R^N)\cap C^{2s+\g}(\R^N)$ if $s<\frac 1 2$, or $u\in \mathcal L^s(\R^N)\cap C^{1, 2s+\g-1}_{loc}(\R^N)\; \textrm{if}\; s\geq \frac 1 2$,  for some $\g>0$. Suppose that $G(u)\in \mathcal L^s(\R^N)$. Then
\begin{equation}\label{eq35}
(-\Delta)^s[G(u)]\leq G'(u)(-\Delta)^s u\quad \textrm{in}\,\;\R^N\,.
\end{equation}
\end{lemma}

\noindent{\it Proof\,.} We can choose, by a suitable convolution, a sequence $\{u_n\}\subset \mathcal S$ uniformly bounded in $C^{2s+\g}_{loc}(\R^N)$, if $s<\frac 1 2$, or in $C^{1, 2s+\g-1}_{loc}(\R^N)\; \textrm{if}\; s\geq \frac 1 2$, for some $\g>0$, with $u_n\to u$ as $n\to\infty$
both in $\mathcal L^s(\R^N)$ and locally uniformly in $\R^N$. Since $G\in C^2(\R;\R)$ and $G(u)\in \mathcal L^s(\R^N)$,
analogously to the proof of \cite[Proposition 2.1.4]{Silv} we have that
$$
(-\Delta)^s u_n\to (-\Delta)^s u, (-\Delta)^s[G(u_n)]\to (-\Delta)^s G(u)\quad \textrm{as}\,\, n\to\infty, 
$$
locally uniformly in $\R^N$. From \cite[Lemma 4.1]{FFV} we have
$$
(-\Delta)^s[G(u_n)]\leq G'(u_n)(-\Delta)^s u_n\quad \textrm{in}\,\;\R^N\,.
$$ 
So, passing to the limit as $n\to \infty$ we get \eqref{eq35}.

\smallskip
\medskip

\begin{proof}[ Proof of Proposition \ref{prop1}] Let $\t\in(0,T)$. Take a nonnegative function 
$$v\in C^2(\bar S_\t)\,\, \text{with}\,\,\operatorname{supp}v(\cdot, t)\,\, \text{compact for each}\,\,t\in [0,\t].$$ 
Moreover, take a function $w\in C(\bar S_\t)\cap  L^1\big((0,\t);\mathcal L^s(\R^N)\big)$ such that for each $t\in(0,\t]$, 
$$\begin{aligned}
&w(\cdot, t)\in \mathcal L^s(\R^N)\cap C^{2s+\g}(\R^N)\,\,\text{if} \,\,s<\frac 1 2,\,\, \text{or}\,\,\\
&w(\cdot, t)\in \mathcal L^s(\R^N)\cap C^{1, 2s+\g-1}_{loc}(\R^N)\; \text{if}\; s\geq \frac 1 2,
\end{aligned}
$$ 
for some $\g>0$. For any $\epsilon\in (0,\tau)$, integrating by parts we have:
\begin{equation}\label{eq36}
\begin{aligned}
& \int_0^\t\int_{\R^N} v \big[-(-\Delta)^s w-  \rho w_t + \left\langle b(x),\nabla w\right\rangle \big]\,  dx dt \\
&= \int_0^\t \int_{\R^N} w\big[- (- \Delta)^s v+  \rho v_t -\left\langle b(x),\nabla v\right\rangle - v\di b\big]\, dx dt\\
&-\int_{\R^N} \rho(x) v(x,\tau)w(x,\tau) dx +\int_{\R^N} \rho(x) v(x,\epsilon) w(x,\epsilon) dx\,.
\end{aligned}
\end{equation}
Let $p\geq 1.$ For any $\a>0$, set
\begin{equation}\label{eq37}
G_\alpha(r):=(r^2+\alpha)^{\frac p 2} \quad \textrm{for all}\;\; r\in\R\,.
\end{equation}
It is easily seen that
\begin{equation}\label{eq38}
\begin{aligned}
&G_\alpha'(r)=pr(r^2+\alpha)^{\frac p2 -1} \\
&G_\alpha''(r)=p(r^2+\alpha)^{\frac p2 -2}[\alpha+r^2(p-1)]\geq 0 \quad \textrm{for all}\;\; r\in\R\,.
\end{aligned}
\end{equation}
By the differential equation in problem \eqref{problema},
\begin{equation}\label{eq39}
\rho [G_\alpha(u)]_t=\rho G_\alpha'(u) u_t=-G_\alpha'(u)[(-\Delta)^s u-\left\langle b,\nabla u\right\rangle ]\quad \textrm{in}\;\, S_T\,.
\end{equation}
From \eqref{eq38}, \eqref{eq39} and Lemma \ref{lemma2} we obtain
\begin{equation}\label{eq310}
\rho [G_\alpha(u)]_t + (-\Delta)^s[G_\alpha(u)]-\left\langle b,\nabla [G_{\alpha}(u)]\right\rangle\le 0\quad \textrm{in}\;\; S_T\,.
\end{equation}
So, from \eqref{eq36} with $w=G_\alpha(u)$ and \eqref{eq310} we obtain
\begin{equation}\label{eq311}
\begin{aligned}
\int_{\R^N} \rho(x) G_\alpha[u(x,\tau)] v(x,\tau) \,dx &\leq  \int_0^\tau\int_{\R^N} G_\alpha(u) \big[- (- \Delta)^s v+\rho v_t-\left\langle b,\nabla v\right\rangle - v\di b \big]\, dx dt \\
& + \int_{\R^N} \rho(x) v(x,\epsilon) G_\alpha[u(x,\epsilon)] \,dx\,.
\end{aligned}
\end{equation}
Letting $\varepsilon\to 0^+$ in \eqref{eq311}, by the dominated convergence theorem,
\begin{equation}\label{eq312}
\begin{aligned}
\int_{\R^N} \rho(x) G_\alpha[u(x,\t)] v(x,\tau)\, dx &\leq  \int_0^\tau \int_{\R^N} G_\alpha(u) \big[- (- \Delta)^s v+\rho  v_t -\left\langle b,\nabla v\right\rangle - v\di b\big] \,dx dt \\ 
&+ \alpha^{p/2} \int_{\R^N} \rho(x) v(x,0)\, dx\,.
\end{aligned}
\end{equation}
Now, letting $\alpha\to 0^+$ in \eqref{eq312}, by the dominated convergence theorem,
\begin{equation}\label{eq313}
\int_{\R^N} \rho(x)  |u(x,\t)|^p v(x,\tau) \, dx \leq  \int_0^\tau \int_{\R^N} |u|^p \big[- (- \Delta)^s v+ \rho  v_t -\left\langle b,\nabla v\right\rangle - v\di b \big] \,dx dt
\end{equation}

For any $R>0$, we can choose
$$
v(x,t):= \phi(x,t)\gamma_R(x)\quad \textrm{for all}\;\; (x, t)\in \bar S_\t\,.
$$
Using the fact that $\phi$ is a supersolution to equation \eqref{eq32c} and $\gamma_R\geq 0$, we obtain
\begin{equation}\label{eq314}
\begin{aligned}
- (- \Delta)^s v&+ \rho v_t -\left\langle b,\nabla v\right\rangle - v\di b\\
&=\gamma_R \left[-(-\Delta)^s \phi + \rho  \phi_t -\left\langle b,\nabla \phi\right\rangle - \phi\di b \right]\\
&\quad\,- \phi(-\Delta)^s\gamma_R + \mathcal B(\phi,\g_R) - \phi\left\langle b, \nabla \gamma_R\right \rangle \\
& \leq - \phi(-\Delta)^s\gamma_R + \mathcal B(\phi,\g_R)- \phi\left\langle b, \nabla \gamma_R\right \rangle  \quad\quad \text{in}\;\; S_\tau\,.
\end{aligned}
\end{equation}
Since $|u|^p\geq 0$, by \eqref{eq313} and \eqref{eq314} we conclude
that
$$
\int_{\R^N} \rho(x)  |u(x,\tau)|^p \phi(x,\tau)\gamma_R(x) dx \leq \int_0^\t \int_{\R^N} |u|^p \big[ - \phi(-\Delta)^s\gamma_R + \mathcal B(\phi,\gamma_R) - \phi\left\langle b, \nabla \gamma_R\right \rangle \big]dx dt.
$$
By using the definitions of $D_+$ and $D_-$ in \eqref{d}, the latter can be rewritten as follows
\begin{equation}\label{equa33}
\begin{aligned}
\int_{\R^N} \rho(x) |u(x,\tau)|^p \phi(x,\tau)\gamma_R(x) dx &\leq \int_0^\t \int_{\R^N} |u|^p \big[ - \phi(-\Delta)^s\gamma_R + \mathcal B(\phi,\gamma_R)\big]\,dxdt \\
&\quad-\int_0^\t \int_{D_+} |u|^p  \phi\left\langle b, \nabla \gamma_R\right \rangle dx dt \\
&\quad -\int_0^\t \int_{D_-} |u|^p  \phi\left\langle b, \nabla \gamma_R\right \rangle dx dt.
\end{aligned}
\end{equation}
Now, observe that, for any $x\in D_-$, for some $\bar C>0$
\begin{equation}\label{eq33c}
\left\langle b(x),\nabla \gamma_R(x)\right\rangle= \gamma'\left(\frac{|x|}{R}\right)\frac{1}{R}\,\left\langle b(x),\frac{x}{|x|}\right\rangle\,\ge -\frac{\bar C}{R}\,\left\langle b(x),\frac{x}{|x|}\right\rangle\ge 0.
\end{equation}
Hence, from \eqref{equa33} and \eqref{eq33c} we have
\begin{equation}\label{eq315}
\begin{aligned}
\int_{\R^N} \rho(x) |u(x,\tau)|^p \phi(x,\tau)\gamma_R(x) dx&\leq \int_0^\t \int_{\R^N} |u|^p \big[ - \phi(-\Delta)^s\gamma_R + \mathcal B(\phi,\gamma_R) \big]dx dt \\
&\quad-\int_0^\t \int_{D_+} |u|^p \phi\left\langle b, \nabla \gamma_R\right \rangle\,dxdt.
\end{aligned}
\end{equation}
Finally, from Lemma \ref{lemma1} with $v=|u|^p$ and the monotone convergence theorem, sending $R\to \infty$ in \eqref{eq315} we get
\begin{equation}\label{eq316}
\int_{\R^N} \rho(x) |u(x,\t)|^p\phi(x,\t)dx\leq 0\,.
\end{equation}
From \eqref{eq316}, \eqref{h3}, \eqref{h0}, since $\phi>0$ in $S_\tau$ and $|u|^p\geq 0$ we infer that $u\equiv 0$ in $S_\tau$. This completes
the proof.
\end{proof}

\medskip

\subsection{Proof of Theorem \ref{teo1}}

Before proving Theorem \ref{teo1}, we need some preliminary results. Observe that the proof of Proposition \ref{prop5} can be found in \cite[Proposition 3.3]{PV1}.

\begin{proposition}\label{prop5}
Let $\tilde w\in C^2([0,\infty))\cap L^\infty((0,\infty)).$ Let
$$w(x):=\tilde w(|x|)\quad \textrm{for all}\;\; x\in \R^N\,.$$
Set $r\equiv |x|$.If
\begin{equation}\label{eq317} 
\tilde w''(r)+ \frac{N-2s+1}{r}\tilde w'(r)\leq 0,
\end{equation}
then $w$ is a supersolution to equation
\begin{equation}\label{eq318}
(-\Delta)^s w\,=\,0\quad \textrm{in}\;\; \R^N\,.
\end{equation}
\end{proposition}

Note that in Proposition \ref{prop5} $w$ is a supersolution to equation \eqref{eq318} in the sense of Definition \ref{defsole} with $c\equiv 0$.
\medskip

In the sequel we shall use the next well-known result, concerning the hypergeometric function $_{2}F_1(a,b,c,s)\equiv F(a,b,c,s)$, with $a,b\in \R, c>0, s\in \R\setminus\{1\}$  (see \cite[Chapters 15.2, 15.4]{DLMF}).

\begin{lemma}\label{lemma3}
The following limits hold true:
\begin{itemize}
\item[(i)] if $c>a+b$, then
\[\lim_{s\to 1^-} F(a,b,c, s)=\frac{\Gamma(c)\Gamma(c-a-b)}{\Gamma(c-a)\Gamma(c-b)}\,;\]
\item[(ii)] if $c=a+b$, then
\[\lim_{s\to 1^-}\frac{F(a,b,c,s)}{-\log(1-s)}=\frac{\Gamma(a+b)}{\Gamma(a)\Gamma(b)}\,;\]
\item[(iii)] if $c<a+b$, then
\[\lim_{s\to 1^-}\frac{F(a,b,c,s)}{(1-s)^{c-a-b}}=\frac{\Gamma(c)\Gamma(a+b-c)}{\Gamma(a)\Gamma(b)}\,.\]
\end{itemize}
\end{lemma}
For further references, observe that
\begin{equation}\label{eq319}
\Gamma(t)>0\quad\textrm{for all}\;\; t>0, \quad \Gamma(t)<0\quad \text{for all}\;\; t\in (-1,0)\,.
\end{equation}
\bigskip

For the proof of Lemma \ref{lemma3}, we refer the reader to \cite[Chapters 15.2, 15.4]{DLMF}.
\bigskip

\begin{proof}[Proof of Theorem \ref{teo1}.] 
Let $\psi=\psi(|x|)$ be defined as in \eqref{eq12}, where $\beta>0$ is a constant to be chosen. Set $r\equiv |x|$. We have:
\begin{align}
&\psi'(r)=-\b r(1+r^2)^{-\left(\frac{\b}2+1\right)}\quad \text{for all}\;\; r>0\,, \label{eq320}\\
&\psi''(r)=\b (1+r^2)^{-\left(\frac{\b}{2}+2\right)}[-1 +(\b+1)r^2]\quad \text{for all}\;\; r>0\,. \label{eq321}
\end{align}
For any~$\lambda>0$ define
$$
\phi(x,t):= e^{-\l t}\psi(r)\quad \text{for all}\;\; (x,t)\in \bar S_T\,.
$$
At first observe that \eqref{eq32d} and \eqref{eq32c} are satisfied. Indeed, one has, for some $C>0$,
$$
\begin{aligned}
\phi(x,t)+|\nabla\phi(x,t)|&\le\, e^{-\lambda t}\left\{\psi(r)+\psi'(r)\right\}\le C\,\psi(r) \quad\quad \text{for all}\,\,(x,t)\in \bar S_T.
\end{aligned}
$$
Furthermore, due to \eqref{h0}, for all $(x,t)\in \bar D_+\times [0,T]$, since $\sigma\le1-\alpha\leq 1$, see \eqref{h3}, one has
$$
\begin{aligned}
\frac{1}{1+|x|}\left\langle b(x),\,\frac{x}{|x|}\right\rangle \phi\,&\le\, e^{-\lambda t}K(1+|x|)^{\sigma-1}\psi(r)\le C\,\psi(r)\,,
\end{aligned}
$$
with $C\geq K.$
\medskip

Now, we want to show that $\phi$ is a supersolution to equation \eqref{eq32b}. To do so, we consider separately the cases $(i)-(iv)$.

\smallskip

Suppose that $(i)$ holds. 
In view of \eqref{eq320} and \eqref{eq321}, we have:
\begin{equation}\label{eq322}
\begin{split}
&\psi''(r)+\frac{N-2s
+1}{r}\psi'(r)\\
&\quad=\b(1+r^2)^{-\left(\frac{\b}2+2\right)}[(\b-N+2s)r^2-(N-2s+2)]\quad \text{for all}\;\; r>0\,.
\end{split}
\end{equation}
Since $0<\b\leq N-2s$, by \eqref{eq322},
\begin{equation}\label{eq323}
\psi''(r)+\frac{N-2s +1}{r}\psi'(r)\leq 0 \quad \text{for all}\;\; r>0\,.
\end{equation}
By Proposition \ref{prop5},
\begin{equation}\label{eq324}
-(-\Delta)^s \phi(x,t)=- e^{-\lambda t}(-\Delta)^s\psi(r)\quad \text{for all}\;\; (x,t)\in \bar S_T\,.
\end{equation}
Furthermore, due to \eqref{h0} and \eqref{eq320}, we get
\begin{equation}\label{eq325}
\begin{aligned}
-\left\langle b(x), \nabla \phi(x,t)\right\rangle & \le e^{-\lambda t}\psi'(r) \left\langle b(x), \frac{x}{|x|}\right\rangle\\
&\le \beta\,K\,e^{-\lambda t}(1+r^2)^{\frac{\sigma}{2}-\frac{\beta}{2}-\frac 12} \quad \quad \text{for all}\,\, (x,t)\in \bar S_T;
\end{aligned}
\end{equation}
and
\begin{equation}\label{eq326}
\begin{aligned}
-\phi(x,t)\di b(x)&\le \phi(x,t)\, [\di b(x)]_{-}\\
&\le e^{-\lambda t}\psi(r)\,K(1+r)^{\sigma-1}\\
&\le K\,e^{-\lambda t}(1+r^2)^{\frac{\sigma}{2}-\frac{\beta}{2}-\frac 12} \quad \quad \quad \text{for all}\,\, (x,t)\in \bar S_T.
\end{aligned}
\end{equation}
From \eqref{eq323}, \eqref{eq324}, \eqref{eq325}, \eqref{eq326} and \eqref{h3}-\eqref{h0}, we obtain
\begin{equation}\label{eq327}
\begin{aligned}
-(-\Delta)^s\phi(x,t) &+\rho(x) \phi_t(x,t)-\left\langle b(x), \nabla \phi(x,t)\right\rangle -\phi(x,t)\di b(x) \\
& \le-  (1+|x|^2)^{-\frac{\beta}{2}} e^{-\lambda t}\left\{\lambda C_0 (1+|x|^2)^{-\frac{\alpha}2}-(\beta K+1)(1+|x|^2)^{\frac{\sigma}{2}-\frac 12}\right\}\\
&\le 0\quad\quad\quad \textrm{for all}\;\; (x,t)\in \bar S_T\,,
\end{aligned}
\end{equation}
provided that 
\[\lambda\geq \frac 1{C_0}(\beta K+1)\,. \]
By \eqref{eq327} and Proposition \ref{prop1}, the conclusion follows, when $(i)$ holds.

\medskip

In order to obtain the thesis of Theorem \ref{teo1} for $\beta>N-2s$ note that (see the proof of Corollary 4.1 in \cite{FerrVerb}) we have:
\begin{equation}\label{eq328}
-(-\Delta)^s\psi(r)= - \check C F(a,b,c, - r^2)\quad \textrm{for all}\;\; r>1\,,
\end{equation}
where $\check C>0$ is a positive constant, and
$$
a=\frac N 2 +s,\quad b= \frac{\b}2+s, \quad c=\frac N2\,.
$$
By Pfaff's transformation,
\begin{equation}\label{eq329}
F(a,b,c, -r^2)=\frac 1{(1+r^2)^b}F\left(c-a, b, c, \frac{r^2}{1+r^2}\right)\quad \textrm{for all}\,\;r>1\,.
\end{equation}

\medskip

 Suppose that $(ii)$ holds. From Lemma \ref{lemma3}-$(i)$, \eqref{eq328} and \eqref{eq329}, for any $\varepsilon>0$, for some $R_\varepsilon>1$, we have:
\begin{equation}\label{eq330}
-(-\Delta)^s\psi(r)\leq \check C(C_1 +\varepsilon)(1+r^2)^{-\left(s+\frac{\beta}2\right)}\quad \text{whenever}\;\;r>R_\varepsilon\,,
\end{equation}
where 
$$
C_1=-\frac{\Gamma(\frac N 2)\Gamma\left(\frac{N-\b}2\right)}{\Gamma\left(\frac{N+s}{2}\right)\Gamma\left(\frac{N-\b}{2}-s\right)}>0
$$ 
(see \eqref{eq319}). Since $\alpha\leq 2s,$ from \eqref{eq330}, \eqref{eq325} and \eqref{eq326}, due to \eqref{h3}-\eqref{h0}, we obtain for all $|x|=r>R_\varepsilon$, $t\in [0,T]$
\begin{equation}\label{eq331}
\begin{aligned}
-(-\Delta)^s \phi(x,t)&+ \rho(x) \phi_t(x,t) -\left\langle b(x), \nabla \phi(x,t)\right\rangle -\phi(x,t)\di b(x) \\
 &\le e^{-\lambda t}(1+|x|^2)^{-\frac{\beta}{2}} \Big\{\check C(C_1+\varepsilon)(1+|x|^2)^{-s}-\lambda C_0(1+|x|^2)^{-\frac{\alpha}2}\\ 
 &\quad\quad+(\beta K+1)(1+|x|^2)^{\frac{\sigma}{2}-\frac 12}\Big\}\\
&\le 0\quad\quad \quad \quad \quad \quad\quad \quad \quad \quad \quad\quad \quad \textrm{for all}\;\; (x,t)\in \bar S_T\,,
\end{aligned}
\end{equation}
provided
\begin{equation}\label{eq332}
\lambda>\frac 2{C_0}\max\{\check C(C_1+\epsilon)\,;\,\beta K+1\}\,.
\end{equation}
On the other hand, for all $|x|\leq R_\varepsilon,\; t\in [0,T]$,
\begin{equation}\label{eq333}
\begin{aligned}
&-(-\Delta)^s \phi(x,t)+ \rho(x) \phi_t(x,t) -\left\langle b(x), \nabla \phi(x,t)\right\rangle -\phi(x,t)\di b(x) \\
& \le-e^{-\lambda t}\left\{-M_{\varepsilon,\beta}+\lambda C_0 (1+R_{\varepsilon}^2)^{-\frac{\beta}{2}-\frac{\alpha}2}-(K\beta+1)\right\}\\
&\le 0
\end{aligned}
\end{equation}
taking
\begin{equation}\label{eq334}
 \lambda>\frac 2{C_0}[M_{\epsilon,\b}+(K\beta+1)](1+R^2_\varepsilon)^{\frac{\beta}2+\frac{\alpha}2},
\end{equation} where
$$
M_{\epsilon,\b}:= \max_{x\in\bar B_{R_\varepsilon}}\big\{\big|-(-\Delta)^s\psi(|x|)\big|\big\}\,.
$$ 
By \eqref{eq331}, \eqref{eq333} the conclusion follows by Proposition \ref{prop1}, when $(ii)$ holds.

\medskip

 Suppose that $(iii)$ holds. From Lemma \ref{lemma3}-$(ii)$ and \eqref{eq329}, for any $\varepsilon>0$, for some $R_\varepsilon>1$, we have:
\begin{equation}\label{eq335}
-(-\Delta)^s\psi(r)\leq \check C(C_2 +\varepsilon)(1+r^2)^{-\left(s+\frac{\b}2\right)}\log(1+r^2)\quad \text{whenever}\;\;|x|>R_\varepsilon\,,
\end{equation}
where
$$C_2=-\frac{\Gamma\left(\frac{\beta}2\right)}{\Gamma(-s)\Gamma\left(\frac{\beta}2+s\right)}>0$$
(see \eqref{eq319}). Since $\alpha<2s$, from \eqref{eq335}, \eqref{eq325}, \eqref{eq326} and \eqref{h3}-\eqref{h0}, we obtain for all $|x|>R_\varepsilon, t\in [0,T]$,
\begin{equation}\label{eq336}
\begin{aligned}
-(-\Delta)^s \phi(x,t)&+ \rho(x) \phi_t(x,t) -\left\langle b(x), \nabla \phi(x,t)\right\rangle -\phi(x,t)\di b(x) \\
&\leq e^{-\lambda t}(1+|x|^2)^{-\frac{\beta}2}\left\{\check C(C_2+\varepsilon)(1+|x|^2)^{-s}\log(1+|x|^2)\right.\\
&\left.\quad\quad-\lambda C_0 (1+|x|^2)^{-\frac{\alpha}2}+(\beta K+1)(1+|x|^2)^{\frac{\sigma}{2}-\frac 12}\right\}\\
&<0\,,
\end{aligned}
\end{equation}
taking a possibly larger $R_\varepsilon>1$, and 
\begin{equation}\label{eq336b}
\lambda>\frac 2{C_0 }\max\left\{\check C(C_2+\varepsilon),\beta K+1\right\}.
\end{equation}
Combining \eqref{eq336} with \eqref{eq333} the conclusion follows, when $(iii)$ holds.

\smallskip

 Finally, suppose that $(iv)$ holds.  From Lemma \ref{lemma3}-$(iii)$ and \eqref{eq329}, for any $\varepsilon>0$, for some $R_\varepsilon>1$, we have:
\begin{equation}\label{eq337}
-(-\Delta)^2\psi(r)\leq \check C(C_3 +\varepsilon)(1+r^2)^{-\left(s+\frac N 2\right)}\quad \text{whenever}\;\;r>R_\varepsilon\,,
\end{equation}
where 
$$C_3=-\frac{\Gamma(\frac N 2)\Gamma\left(\frac{\beta-N}2\right)}{\Gamma(-s)\Gamma\left(\frac{\beta}{2}+s\right)}>0
$$
(see \eqref{eq319}). Since $\alpha+\beta\leq 2s+N$, from \eqref{eq337}, \eqref{eq325}, \eqref{eq326}, and due to \eqref{h3}-\eqref{h0}, we obtain for all $|x|>R_\varepsilon$, $t\in [0,T]$
\begin{equation}\label{eq338}
\begin{aligned}
-(-\Delta)^s \phi(x,t)&+ \rho(x) \phi_t(x,t) -\left\langle b(x), \nabla \phi(x,t)\right\rangle -\phi(x,t)\di b(x) \\
&\leq e^{-\lambda t}(1+|x|^2)^{-\frac{\beta}2}\left\{\check C(C_3+\varepsilon)(1+|x|^2)^{-s-\frac{N}{2}+\frac{\beta}2}\right.\\
&\left.\quad\quad-\lambda C_0 (1+|x|^2)^{-\frac{\alpha}2}+(\beta K+1)(1+|x|^2)^{\frac{\sigma}{2}-\frac 12}\right\}\\
&<0\,,
\end{aligned}
\end{equation}
provided
\begin{equation}\label{eq339}
\lambda>\frac 2{C_0}\max\{\check C(C_3+\epsilon)\,;\,\beta K+1\}\,.
\end{equation}
On the other hand, \eqref{eq333} holds true, provided \eqref{eq334} is satisfied. In view of \eqref{eq333} and \eqref{eq338}, the conclusion follows by Proposition \ref{prop1}, when $(iv)$ holds. This completes the proof.
\end{proof}

\subsection{Proof of Proposition \ref{cor2}}

To show nonuniqueness for problem \eqref{problema} with $b$ as in \eqref{h2}, we use the following result.

\begin{proposition}\label{prop32} Let $\rho, b\in C^1(\mathbb R^N), \rho>0$. 
 If there exists a viscosity supersolution of problem 
\begin{equation}
\begin{aligned}\label{eq26}
-(-\Delta)^s V+\left\langle b,\,\nabla V\right\rangle&=-\rho \quad\quad \text{in}\,\,\, \R^N,\\
V>0 \quad \text{ in }\,\, \mathbb R^N, \quad \lim_{|x|\to+\infty}V(x)&=0\,,
\end{aligned}
\end{equation}
then there exist infinitely many bounded solutions $u$ of problem \eqref{problema}. In particular, for any $g\in C([0,T])$, $g(0)=0$, there exists a solution $u$ to problem \eqref{problema} such that
$$
\lim_{|x|\to\infty} u(x,t)=g(t) \quad \text{uniformly with respect to}\,\,\,t\in[0,T].
$$
\end{proposition}
Proposition \ref{prop32} can be proved by minor changes in the proof of \cite[Theorem 2.7]{PV2}. In particular, in view of  the hypotheses on $\rho$ and $b$, by regularity results, the  constructed solutions satisfy problem \eqref{problema} in the sense of Definition \ref{defsolp} (see e.g. \cite{Silv,Silv2,Silv3}).  
\medskip

Due to Proposition \ref{prop32}, in order to prove Proposition \ref{cor2}, it is enough to show that such a supersolution $V$ to \eqref{eq26} exists, provided that assumption \eqref{h2} is satisfied. 

\begin{lemma}\label{lemma5} Let $\rho, b\in C(\mathbb R^N)$, $\rho>0$ in $\R^N$. Assume that \eqref{h2} holds. Then there exists a viscosity supersolution $V>0$ of problem \eqref{eq26}.

\end{lemma}

\begin{proof} For any $C>0$ and $\beta>0$ we define the function
\begin{equation}\label{eq340}
V_1(x):= C|x|^{-\beta}\quad \text{for any}\,\,x\in\R^N\setminus\{0\},
\end{equation}
where $C>0$ and $\beta$ have to be chosen. 
Note that (see Section \ref{sec1})
$$
(-\Delta)^sV_1(x)=C_{N,s}\,\int_{\R^N}\frac{V_1(x)-V_1(y)}{|x-y|^{N+2s}}\,dy\quad \text{for all}\,\,x\in\R^N\setminus\{0\}.
$$
Moreover we have
$$
\hat{V_1}(\xi)=C_{\beta}C|\xi|^{-N+\beta}
$$
for some $C_{\beta}>0$. Hence, by \eqref{eq23}
$$
(-\Delta)^sV_1(x)=CC_{\beta}\left(\mathcal{F}^{-1}|\xi|^{2s-N+\beta}\right)(x)=CC_{\beta}C_{\beta+2s}|x|^{-\beta-2s}\quad \text{for all}\,\,x\in\R^N\setminus\{0\}.
$$
Thus, in view of \eqref{h2}, we have
\begin{equation}\label{eq342}
\begin{aligned}
-(-\Delta)^s V_1(x)+\left\langle b(x),\,\nabla V_1(x)\right\rangle\,&\le\, -C C_{\beta}C_{\beta+2s}|x|^{-\beta-2s}-\beta C |x|^{\beta-1}\left\langle b(x),\,\frac{x}{|x|}\right\rangle\\
&\le\,-C C_{\beta}C_{\alpha+2s}|x|^{-\beta-2s}-\beta C K |x|^{\sigma-1-\beta} \\
&\le -\beta C K |x|^{\sigma-1-\beta},
\end{aligned}
\end{equation}
for all $x\in\R^N\setminus B_{R_0}$, with $R_0$ as in \eqref{h2}. Now, from \eqref{eq342} it follows that 
\[-(-\Delta)^s V_1(x)+\left\langle b(x),\,\nabla V_1(x)\right\rangle\,\le\,-\rho \quad \text{for all}\,\,\,x\in\R^N\setminus\bar{B}_{R_0};\]
provided that
$$
\beta=\sigma-1-\alpha>0, 
$$
and $C>0$ is sufficiently large.

Define 
\[\hat V(x):= C_2 (R_0^2-|x|^2)_+^{\frac s2}, \quad x\in \mathbb R^N\,.\]
From the results in \cite{Geetor} it follows that, for a suitable $C_2>0$, $\hat V$ solves 
\begin{equation}\label{eq210}
\begin{cases}
-(-\Delta)^s \hat V = -1 & \text{ in } B_{R_0}\\
\hat V =0 & \text{ in } \mathbb R^N\setminus B_{R_0}\,.
\end{cases}
\end{equation}
In view of \eqref{h2}, from \eqref{eq210} it immediately follows that, for $C_2\geq \max_{\bar B_{R_0}}\rho$,
 \[V_2:= C_2 \hat V\]
 solves 
\[\begin{cases}
-(-\Delta)^s V_2 + \langle b(x), \nabla V_2(x) \rangle\leq  -\rho & \text{ in } B_{R_0}\\
V_2 =0 & \text{ in } \mathbb R^N\setminus B_{R_0}\,.
\end{cases}
\]
By using $V_1$ and $V_2$, by the same arguments as in the proof of \cite[Proposition 3.2]{PV2} we get the thesis.
\end{proof}

\begin{proof}[Proof of Proposition \ref{cor2}]
By simply combining Proposition \ref{prop32} and Lemma \ref{lemma5}, the result follows.
\end{proof}

\section{Elliptic equations: proofs}\label{pe}\setcounter{equation}{0}

Firstly we state a general criterion for uniqueness of nonnegative solutions to equation \eqref{elliptic} in $L^1_{\zeta}(\R^N)$. We suppose that there exists a positive function $\zeta\in C^2(\R^N)$, which solves, for some $p>1$
\begin{equation}\label{eq40}
-(-\Delta)^s\zeta - \left \langle b,\,\nabla \zeta \right \rangle -\zeta \di b\, - \rho\,p\,c \,\zeta < 0 \quad \textrm{in}\;\; \R^N\,.
\end{equation}
Such inequality is meant in the sense that in Definition \ref{defsole}-$(ii)$, instead of $``="$ we have $``<"$.

\begin{proposition} \label{prop3}
Let assumptions \eqref{h3}, \eqref{h0} and \eqref{h1} be satisfied. Let $u$ be a solution to equation \eqref{elliptic} with $|u|^p\in \mathcal L^s(\R^N),$ for some $p\geq 1$. Assume that there exists a positive function $\zeta\in C^2(\R^N)$, which solves \eqref{eq40}, and satisfies
\begin{equation}\label{eq42d}
\zeta(x)+|\nabla \zeta(x)|\leq C \psi(x)\quad \textrm{for all}\;\; x\in \R^N;
\end{equation}
and
\begin{equation}\label{eq42c}
\frac{1}{1+|x|}\left\langle b(x),\frac{x}{|x|}\right\rangle\zeta(x)\leq C \psi(x)\quad \textrm{for all}\;\; x\in D_+;
\end{equation}
for some constant $C>0$ and for $\psi$ as in \eqref{eq12}. If $u\in L^p_\psi(\R^N)$, then
$$
u\equiv 0\quad \textrm{in}\;\; \R^N\,.
$$
\end{proposition}

\subsection{Proof of Proposition \ref{prop3}}
Analogously to Lemma \ref{lemma1} the next lemma can be shown.
\begin{lemma}\label{lemma4}
Assume \eqref{h0} and \eqref{h1}. Let $\zeta\in C^2(\R^N)$, $\zeta>0$; suppose that \eqref{eq42d} and \eqref{eq42c} are satisfied. Let $v\in L^1_{\psi}(\R^N)$. Then
\begin{equation}\label{eq41bis}
\int_{\R^N}|v(x)|\zeta(x)|(-\Delta)^s\gamma_R(x)|\,dx + \int_{\R^N}|v(x)|\,|\mathcal B(\zeta, \gamma_R)(x)|dx \to 0, 
\end{equation}
as $R\to \infty\,,$ and
\begin{equation}\label{eq41tris}
\lim_{R\to\infty}\int_{D_+}|v(x)|\zeta(x)\left\langle b(x)\,,\,\nabla\gamma_R(x)\right\rangle\,dx = 0, 
\end{equation}
where $D_+$ is defined in \eqref{d}.
\end{lemma}

\begin{remark}\label{rem31}
Observe that the quantity into the brackets of formula \eqref{eq41tris} is negative (see Remark \ref{rem21} for the explanation).
\end{remark}

\begin{proof}[ Proof of Proposition \ref{prop3}]\,. Take a function $v\in C^2(\R^N)$ with $\operatorname{supp} v$ compact. Moreover, take a function $w\in C(\R^N)$ such that $w\in \mathcal L^s(\R^N)\cap C^{2s+\g}_{loc}(\R^N)$ if $s<\frac 1 2$, or  $w\in \mathcal L^s(\R^N)\cap C^{1, 2s+\g-1}_{loc}(\R^N)$ if $s\geq\frac 1 2$, for some $\gamma>0$. Integrating by parts we have:
\begin{equation}\label{eq41}
\begin{aligned}
\int_{\R^N} &v \left[-(-\Delta)^s w +\left\langle b\,,\,\nabla w\right\rangle - \rho(x)  c(x) w \right]\,  dx  \\
&= \int_{\R^N} w\left[- (- \Delta)^s v - \left\langle b\,,\,\nabla v\right\rangle - v\di b- \rho(x)  c(x) v\right]\, dx.
\end{aligned}
\end{equation}
Let $G_\alpha$ be defined as in \eqref{eq37}. Due to \eqref{elliptic} and \eqref{eq39} we obtain
\begin{equation}\label{eq42}
\begin{aligned}
(-\Delta)^s[G_\alpha(u)] &- \left\langle b\,,\,\nabla [G_{\alpha}(u)] \right\rangle + \rho  c \,G_\alpha(u) \\
& \le G_{\alpha}'(u)(-\Delta)^su-G_{\alpha}'(u)\left\langle b\,,\,\nabla u \right\rangle + \rho c G_\alpha(u) \\
& = p\,u (u^2+\alpha)^{\frac p2-1}(-\Delta)^s u - p\,u (u^2+\alpha)^{\frac p2-1}\left\langle b\,,\,\nabla u \right\rangle  + \rho\, c\, G_{\alpha}(u) \\
&\quad\quad+  \rho\,  c\,p\,u^2(u^2+\alpha)^{\frac p2-1} - \rho\, c\,p\,u^2(u^2+\alpha)^{\frac p2-1}\\
&=p\,u (u^2+\alpha)^{\frac p2 -1}[(-\Delta)^s u -\langle b\,,\,\nabla u\rangle + \rho\,  c\, u] \\
&\quad\quad +  \rho\, c(u^2+\a)^{\frac p2-1}\left[u^2+\alpha-p\,u^2\right]\\
&=\rho\, c(u^2+\a)^{\frac p2-1}\left[u^2(1-p)+\alpha\right]\quad \quad\quad\textrm{in}\;\; \R^N\,.
\end{aligned}
\end{equation}
From \eqref{eq41} with $w=G_\alpha(u)$ and \eqref{eq42} it follows that
\begin{equation}\label{eq43}
\begin{aligned}
\int_{\R^N} v (u^2+\alpha)^{\frac p2-1} &\rho(x) c(x)[(p-1)u^2- \alpha]\\
&\leq \int_{\R^N} G_\alpha(u)\big[- (- \Delta)^s v- \left\langle b\,,\,\nabla v\right\rangle - v\di b-\rho(x) c\, v\big]\, dx.
\end{aligned}
\end{equation}
Letting $\a\to 0^+$ in \eqref{eq43}, by the dominated convergence
theorem we get
\begin{equation}\label{eq44}
\int_{\R^N} |u|^p\big[(- \Delta)^s v+\left\langle b\,,\,\nabla v\right\rangle + v\di b+p\,\rho\, c\, v\big]\, dx\,\le \,0.
\end{equation}
For any $R>0$, we can choose
\[v(x):= \zeta(x)\g_R(x)\quad \textrm{for all}\;\; x\in \R^N\,,\]
where $\zeta$ satisfies assumption \eqref{eq40}. Then, we obtain
\begin{equation}\label{eq45}
\begin{aligned}
(-\Delta)^s v+\left\langle b\,,\,\nabla v\right\rangle &+ v\di b+\,p\,\rho\, c v\\
&=\gamma_R \left[(-\Delta)^s \zeta +\left\langle b\,,\,\nabla \zeta\right\rangle + \zeta\di b + p\,\rho\, c \,\zeta\right]\\
&\quad+ \zeta(-\Delta)^s\gamma_R -\mathcal B(\zeta,\gamma_R) +\zeta\,\langle b\,,\,\nabla \gamma_R\rangle \quad \quad\text{in}\;\; \R^N\,.
\end{aligned}
\end{equation}
By arguing as in the proof of Proposition \ref{prop1}, due to \eqref{eq44}, \eqref{eq45} and the definitions of $D_+$ and $D_-$ in \eqref{d}, we obtain
\begin{equation}\label{eq46}
\begin{aligned}
\int_{\R^N} |u(x)|^p \gamma_R& \left[(-\Delta)^s \zeta +\left\langle b\,,\,\nabla \zeta\right\rangle + \zeta\di b + p\, \rho\, c \,\zeta\right]\, dx \\
&\le \int_{\R^N} |u|^p \left[-\zeta(-\Delta)^s\gamma_R +\mathcal B(\zeta,\gamma_R) -\zeta\,\langle b\,,\,\nabla \gamma_R\rangle\right]\,dx\\
&\le \int_{\R^N} |u|^p \left[-\zeta(-\Delta)^s\gamma_R +\mathcal B(\zeta,\gamma_R)\right] \,dx\\
&\quad-\int_{D_+}|u|^p\zeta\,\langle b\,,\,\nabla \gamma_R\rangle\,dx -\int_{D_-}|u|^p\zeta\,\langle b\,,\,\nabla \gamma_R\rangle\,dx\,.
\end{aligned}
\end{equation}
Observing that, for all $x\in D_-$
$$
- \langle b(x)\,,\,\nabla \gamma_R(x)\rangle \,\le 0,
$$
we rewrite \eqref{eq46} as follows
\begin{equation}\label{eq46a}
\begin{aligned}
\int_{\R^N} |u|^p \gamma_R& \left[(-\Delta)^s \zeta +\left\langle b\,,\,\nabla \zeta\right\rangle + \zeta\di b + p\,\rho\, c \,\zeta\right]\, dx \\
&\le \int_{\R^N} |u|^p \left[-\zeta(-\Delta)^s\gamma_R +\mathcal B(\zeta,\gamma_R)\right] \,dx\\
&\quad -\int_{D_+}|u|^p\zeta\,\langle b\,,\,\nabla \gamma_R\rangle\,dx\,.
\end{aligned}
\end{equation}
Hence, due to \eqref{eq42d} and \eqref{eq42c}, we can apply Lemma \ref{lemma4}. From Lemma \ref{lemma4} and the monotone convergence theorem, sending $R\to \infty$ in \eqref{eq46a} we get
\begin{equation}\label{eq47}
\int_{\R^N} |u|^p \left[(-\Delta)^s \zeta +\left\langle b\,,\,\nabla \zeta\right\rangle + \zeta\di b + p\, \rho\, c \,\zeta\right]\,dx\leq 0\,.
\end{equation}
From \eqref{eq47} and \eqref{eq40}, since $|u|^p\geq 0$, we can infer that $u\equiv 0$ in $\R^N$. This completes the proof.
\end{proof}

\medskip

\subsection{Proof of Theorem \ref{teo3}}

\begin{proof}[Proof of Theorem \ref{teo3}\,.] 
Let $\psi$ be defined by \eqref{eq12}. From the same arguments as in the proof of Theorem \ref{teo1} we can infer that $\psi$ fulfills \eqref{eq42d} and \eqref{eq42c}, for $\beta>0$. Moreover, $\psi$ solves \eqref{eq40}, for properly chosen $\beta>0$. Note that to do this, we require that $p c_0$, where $c_0$ is defined in \eqref{eq27}, satisfies the same conditions as $\lambda$ in the proof of Theorem \ref{teo1}. To be specific, let us consider all the cases $(i)-(iv)$ as well as in the proof of Theorem \ref{teo1}.
\begin{itemize}
\item Let $(i)$ holds. We should require 
 $$p\, c_0>\frac 1{C_0}(\beta K+1),$$
for $K$ as in \eqref{h0}. 

\item Let $(ii)$ hold. We require that, for some $\varepsilon>0$,
\begin{equation}\label{eq48}
p \, c_0  >\frac2{C_0}\max\left\{\check C(C_1+\epsilon),\, \beta K+1,\, M_{\varepsilon, \beta}(1+R^2_\varepsilon)^{\frac{\beta}2+\frac{\alpha}2}\right\}\,,
\end{equation}
(see \eqref{eq332}, \eqref{eq334}).
\item Let $(iii)$ hold. It is sufficiently to have, for some $\varepsilon>0$
\begin{equation}\label{eq49}
p\, c_0  >\frac 2{C_0}\max\left\{\check C(C_2+\varepsilon),\, \beta K+1,\, M_{\varepsilon,\beta}(1+R^2_\varepsilon)^{\frac{\b}2+\frac{\alpha}{2}}\right\}\,,
\end{equation}
(see \eqref{eq336b}, \eqref{eq334}). 
\item Finally, Let $(iv)$ hold. We ask that, for some $\varepsilon>0$
\begin{equation} \label{eq410} 
p\, c_0 >\frac 2{C_0}\max\left\{\check C(C_3+\varepsilon),\, \beta K+1,\, M_{\varepsilon,\beta}(1+R^2_\varepsilon)^{\frac{\b}2+\frac{\alpha}{2}}\right\}\,,
\end{equation}
(see \eqref{eq339}, \eqref{eq334}).
\end{itemize}
Thus, by Proposition \ref{prop3} the conclusion follows. 
\end{proof}

\subsection{Proof of Proposition \ref{cor4}}
To show nonuniqueness for problem \eqref{elliptic} with $b$ as in \eqref{h2}, we use the following result.

\begin{proposition}\label{prop42}  Let $\rho$, $b$, $c\in C^1(\mathbb R^N)$, $\rho>0$, $c\geq 0$ in $\R^N$. 
Suppose that there exists a viscosity supersolution $h$ of equation 
\begin{equation}\label{eq100}
-(-\Delta)^s h+\left\langle b,\,\nabla h\right\rangle  =-\rho \quad\quad \text{in}\,\,\, \R^N\,.
\end{equation}
such that 
\[h>0 \quad \text{ in }\,\, \mathbb R^N, \quad \lim_{|x|\to +\infty} h(x)=0\,.\]
Then there exist infinitely many bounded solutions $u$ of problem \eqref{elliptic}. In particular, for any $\gamma\in \mathbb R$, there exists a solution $u$ to equation \eqref{elliptic} such that
$$
\lim_{|x|\to\infty} u(x)=\gamma.
$$
\end{proposition}
Proposition \ref{prop42} can be proved by minor changes in the proof of \cite[Theorem 2.10]{PV2}.

\begin{proof}[Proof of Proposition \ref{cor4}] By Lemma \ref{lemma5} and Proposition \ref{prop42}, the conclusion follows. 
\end{proof}

\end{document}